\documentstyle[amssymb,12pt]{amsart}%\input amssym.def
\newtheorem{prop}{Proposition}[section]
\newtheorem{defin}{Definition}[section]
\newtheorem{lemma}{Lemma}[section]
\newtheorem{thm}{Theorem}[section]

\theoremstyle{remark}
\newtheorem{remark}{Remark}

\begin{document}

\newcommand{\nc}{\newcommand} \nc{\on}{\operatorname}

\nc{\cA}{{\cal A}}\nc{\cB}{{\cal B}}\nc{\cC}{{\cal C}}
\nc{\cE}{{\cal E}}\nc{\cG}{{\cal G}}\nc{\cH}{{\cal H}}
\nc{\cX}{{\cal X}}\nc{\cR}{{\cal R}}\nc{\cL}{{\cal L}}
\nc{\cK}{{\cal K}}\nc{\cO}{{\cal O}} \nc{\cF}{{\cal F}}
\nc{\cS}{{\cal S}}

\nc{\sh}{\on{sh}}\nc{\Id}{\on{Id}}\nc{\Diff}{\on{Diff}}
\nc{\Sym}{\on{Sym}}\nc{\ad}{\on{ad}}\nc{\Der}{\on{Der}}
\nc{\End}{\on{End}}\nc{\res}{\on{res}}
\nc{\Imm}{\on{Im}}\nc{\limm}{\on{lim}}\nc{\Ad}{\on{Ad}}
\nc{\Hol}{\on{Hol}}\nc{\Det}{\on{Det}}\nc{\Hom}{\on{Hom}}
\nc{\dimm}{\on{dim}}\nc{\Ker}{\on{Ker}}\nc{\Coeff}{\on{Coeff}}
\nc{\Rep}{\on{Rep}}

\nc{\al}{\alpha}\nc{\de}{\delta}\nc{\si}{\sigma}\nc{\ve}{\varepsilon}
\nc{\la}{{\lambda}}\nc{\g}{\gamma}\nc{\eps}{\epsilon}

\nc{\La}{\Lambda}

\nc{\pa}{\partial}

\nc{\CC}{{\mathbb C}}\nc{\ZZ}{{\mathbb Z}}\nc{\NN}{{\mathbb
N}}\nc{\FF}{{\mathbb F}} \nc{\AAA}{{\mathbb A}}

\nc{\G}{{\mathfrak g}}\nc{\A}{{\mathfrak a}}
\nc{\C}{{\mathfrak c}} \nc{\HH}{{\mathfrak h}}
\nc{\N}{{\mathfrak n}}\nc{\B}{{\mathfrak b}}
\nc{\SL}{{\mathfrak{sl}}}

\nc{\wt}{\widetilde} \nc{\wh}{\widehat}

\nc{\ux}{{\underline{x}}}\nc{\uy}{{\underline{y}}}

\nc{\bn}{\begin{equation}}\nc{\en}{\end{equation}}

% ****** GISPIC **********
%
%** by GISLI MASON *******
%
%**for commutative diagrams
%

\newcommand{\ldar}[1]{\begin{picture}(10,50)(-5,-25)
\put(0,25){\vector(0,-1){50}}
\put(5,0){\mbox{$#1$}} 
\end{picture}}

\newcommand{\lrar}[1]{\begin{picture}(50,10)(-25,-5)
\put(-25,0){\vector(1,0){50}}
\put(0,5){\makebox(0,0)[b]{\mbox{$#1$}}}
\end{picture}}

\newcommand{\luar}[1]{\begin{picture}(10,50)(-5,-25)
\put(0,-25){\vector(0,1){50}}
\put(5,0){\mbox{$#1$}}
\end{picture}}

\title[Quantum homogeneous spaces and quasi-Hopf algebras]
{Quantum homogeneous spaces \\ and quasi-Hopf algebras}

\author{B. Enriquez}

\address{B.E.: D\'epartement de Math\'ematiques et Applications (UMR
8553 du CNRS), Ecole Normale Sup\'erieure, 45 rue d'Ulm, 75005 Paris,
France}

\author{Y. Kosmann-Schwarzbach}

\address{Y.K.-S.: Centre de Math\'ematiques (UMR 7640 du CNRS), 
Ecole Polytechnique, 91128 Palaiseau, France}

\dedicatory{To the memory of Mosh\'e Flato}

\date{January 2000}

\begin{abstract}
We propose a formulation of the quantization problem of Manin quadruples, 
and show that a solution to this problem yields a quantization of the 
corresponding Poisson homogeneous spaces. We then solve both quantization 
problems in an example related to quantum spheres.
\end{abstract}

\maketitle

\section*{Introduction}

According to a theorem of Drinfeld, formal Poisson homogeneous spaces 
over a formal Poisson-Lie group $G_+$ with Lie algebra $\G_+$ correspond
bijectively   to $G_+$-conjugation classes of Lagrangian (i.e., maximal
isotropic) Lie subalgebras $\HH$ of the double Lie algebra $\G$ of $\G_+$. 
The  formal Poisson homogenous space is then $G_+ / (G_+ \cap H)$, 
where $H$ is the formal Lie group with Lie algebra $\HH$. The
corresponding quantization problem is to deform the algebra of functions
over  the  homogeneous space to an algebra-module over the quantized
enveloping algebra of $\G_+$. 

In this paper, we show that there is a Poisson homogeneous structure
on the formal homogeneous space $G/H$, such that the embedding of $G_+
/ (G_+ \cap H)$ in $G/H$ is Poisson, where $G$ is the formal Lie group
with Lie algebra $\G$. It is therefore natural to seek a quantization
of the function algebra of $G_+ / (G_+ \cap H)$ as a quotient of a
quantization of $G/H$.

The data of $(\G,\HH)$ and the $r$-matrix of $\G$ constitute an
example of a quasitriangular Manin pair (see Section
\ref{terminology}). We introduce the notion of the quantization of a
quasitriangular Manin pair, which consists of a quasitriangular Hopf
algebra quantizing the Lie bialgebra $\G$, quasi-Hopf algebras
quantizing the Manin pair $(\G,\HH)$, and a twist element relating
both structures.  We then show (Theorem \ref{thm:main}) that this data
gives rise to a quantization of $G/H$.

The quadruple $(\G,\G_+,\G_-,\HH)$ formed by adjoining $\HH$ to the
Manin triple of $\G_+$ is called a Manin quadruple.  The quantization
of a Manin quadruple is the additional data of Hopf algebras
quantizing the Manin triple $(\G,\G_+,\G_-)$, subject to compatibility
conditions with the quantization of the underlying quasitriangular
Manin pair.  We show that any quantization of a given quadruple gives
rise to a quantization of the corresponding homogeneous space (Theorem
\ref{thm:bis}).

Finally, in Section \ref{sect:spheres}, we explicitly solve the 
problem  of quantizing a Manin quadruple in a situation related to 
quantum  spheres.

\section{Manin quadruples and Poisson homogeneous spaces} 

In this section, we define Lie-algebraic structures, i.e.,
quasitriangular Manin pairs and Manin quadruples, and the
Poisson homogeneous spaces naturally associated to them.

\subsection{Quasitriangular Manin pairs} \label{terminology}

Recall that a {\it quasitriangular Lie bialgebra} is a triple $(\G,r,\langle\
,\ \rangle_{\G})$, where $\G$ is a complex Lie algebra, $\langle \ ,\
\rangle_\G$ is a nondegenerate invariant symmetric  bilinear form on
$\G$, and $r$ is an element of $\G\otimes \G$ such that $r + r^{(21)}$
is the symmetric element of $\G\otimes\G$ defined by $\langle\ ,\
\rangle_\G$, and $r$ satisfies the classical Yang-Baxter equation,
$[r^{(12)},r^{(13)}] + [r^{(12)}, r^{(23)}] + [r^{(13)},r^{(23)}] = 0$.
(Such Lie bialgebras are also  called factorizable.)
 
Assume that $\G$ is finite-dimensional, and let $G$ be a Lie group with
Lie algebra $\G$. Then $G$ is equipped with the Poisson-Lie bivector 
$P_G = r^{L} - r^{R}$, where, for any element $a$ of $\G\otimes\G$, we
denote the right- and left-invariant $2$-tensors on $G$ corresponding to
$a$ by $a^L$ and $a^R$. If $\G$ is an arbitrary Lie algebra, the same
statement holds for its formal Lie group.

We will call the pair $(\G,r,\langle\ ,\ \rangle_{\G},\HH)$ of a
quasitriangular Lie bialgebra $(\G,r,\langle\ ,\ \rangle_{\G})$ and a
Lagrangian Lie subalgebra $\HH$ in $\G$ a {\it quasitriangular Manin
pair}.

Assume that $(\G,r,\langle\ ,\ \rangle_{\G},\HH)$ is a quasitriangular
Manin pair, and that $L$ is a Lagrangian complement of $\HH$ in $\G$.
Let $(\eps^i)$ and $(\eps_i)$ be dual bases of $\HH$ and $L$, and set 
$r_{\HH,L} = \sum_i \eps^i\otimes \eps_i$. The restriction to $L$
of the Lie bracket of $\G$ followed by the projection to the first factor in 
$\G = \HH \oplus L$ yields an element $\varphi_{\HH,L}$ of $\wedge^3\HH$. 
Let us set $f_{\HH,L} = r_{\HH,L} - r$; then $f_{\HH,L}$
belongs to $\wedge^2\G$.  

Then the twist of the Lie bialgebra $(\G,\pa r)$ by $f_{\HH,L}$ is the 
quasitriangular Lie quasi-bialgebra $(\G,\pa r_{\HH,L},
\varphi_{\HH,L})$.  The cocycle $\pa r_{\HH,L}$ maps $\HH$ to
$\wedge^2\HH$, so $(\HH,(\pa r_{\HH,L})_{|\HH},\varphi_{\HH,L})$ is a 
sub-Lie quasi-bialgebra of $(\G,\pa r_{\HH,L},\varphi_{\HH,L})$.

\subsection{Manin quadruples} 

\subsubsection{Definition}

A {\it Manin quadruple} is a quadruple $(\G,\G_+,\G_-,\HH)$, where $\G$ is a
complex Lie algebra, equipped with a nondegenerate invariant symmetric 
bilinear form $\langle \ ,\ \rangle_\G$, and $\G_+,\G_-$ and $\HH$ are
three Lagrangian subalgebras of $\G$, such that $\G$ is equal to the
direct sum $\G_+ \oplus\G_-$ (see \cite{EK}).

In particular, $(\G,\G_+,\G_-)$ is a Manin triple. This implies that for
any Lie group $G$ with subgroups $G_+,G_-$ integrating $\G,\G_+,\G_-$,
we obtain a Poisson-Lie group structure $P_G$ on $G$, such that $G_+$
and $G_-$ are Poisson-Lie subgroups of $(G,P_G)$ (see \cite{Dr:QG}). We
denote the corresponding Poisson structures on $G_+$ and $G_-$ by
$P_{G_+}$ and $P_{G_-}$.

Any Manin quadruple $(\G,\G_+,\G_-,\HH)$ gives rise to a quasitriangular
Manin pair $(\G,r,\langle\ ,\ \rangle_\G,\HH)$, where  $r =
r_{\G_+,\G_-}= \sum_i e^i\otimes e_i$, and $(e^i),(e_i)$ are dual bases
of $\G_+$ and $\G_-$.

\subsubsection{Examples}

In the case where $(\G,\G_+,\G_-)$ is the Manin triple associated with
the Sklyanin structure on a semisimple Lie group $G_+$, the Manin
quadruples were classified in \cite{Karol}. (See also \cite{Delorme}.)
It was shown in \cite{Lu2} that the Lagrangian subalgebras  $\HH$ of
$\G$ such that the intersection $\G_+\cap \HH$ is a Cartan  subalgebra
of $\G_+$ correspond bijectively to the  classical dynamical
$r$-matrices for $\G_+$. We will treat the  quantization of this example
in Section \ref{sect:examples}, in  the case where $\G_+ = \SL_2$.  

In \cite{ER}, the following class of Manin quadruples was studied. Let
$\bar\G$ be a semisimple Lie algebra with nondegenerate invariant
symmetric bilinear form $\langle\ ,\ \rangle_{\bar\G}$, and Cartan 
decomposition $\bar\G =\bar\N_+\oplus\bar\HH\oplus \bar\N_-$.  Let $\cK$
be a commutative ring equipped with a nondegenerate symmetric  bilinear
form $\langle\ ,\ \rangle_\cK$. Assume that $R\subset \cK$ is a
Lagrangian subring of $\cK$, with Lagrangian complement $\La$.  Let us
set  $\G = \bar\G\otimes\cK$, let us equip $\G$ with the bilinear form
$\langle\ ,\ \rangle_{\bar\G} \otimes \langle\ ,\ \rangle_\cK$ and let
us set 
\begin{equation} \label{infinite:quadruple}
 \G_+ = (\bar\HH \otimes R) \oplus (\bar\N_+ \otimes \cK) , \quad \G_- =
 (\bar\HH \otimes \La) \oplus (\bar\N_- \otimes \cK), \quad \HH =
 \bar\G\otimes R.   
\end{equation} 
More generally, extensions of these Lie algebras, connected with the
additional data of a derivation $\pa$ of $\cK$ leaving $\langle\ ,\
\rangle_\cK$ invariant and preserving $R$, were considered in \cite{ER}. 
Examples of quadruples $(\cK,\pa,R,\langle\ ,\ \rangle_{\cK})$, where
$\cK$ is an infinite-dimensional vector space, arise in the
theory of complex curves.

\subsection{Formal Poisson homogeneous spaces}

In what follows, all homogeneous spaces will be formal, so if $\A$ is a
Lie algebra and $A$ is the associated formal group, the function ring of
$A$ is $\cO_A := (U\A)^*$, and if $\B$ is a Lie subalgebra of $\A$,
and $B$  is the associated formal group, the function ring of $A/B$ is 
$\cO_{A/B} :=\left( U\A / (U\A) \B \right)^*$. 

We will need the following result on formal homogeneous spaces. 

\begin{lemma} \label{lemma:surj}
Let $\A$ be a Lie algebra, let $\A_+$ and $\B$ be Lie subalgebras of
$\A$,  and let $A,A_+$ and $B$ be the associated formal groups.  The
restriction map $(U\A / (U\A)\B)^* \to (U\A_+ / (U\A_+(\A_+\cap  \B))
)^*$ is a surjective morphism of  algebras from  $\cO_{A/B}$ to
$\cO_{A_+/(A_+\cap B)}$. 
\end{lemma}

{\em Proof.} Let $L_+$ be a complement of $\A_+\cap \B$ in 
$\A_+$, and let $L$ be a complement of $\B$ in $\A$, containing 
$L_+$. When $V$ is a vector space, we denote by $S(V)$ its symmetric 
algebra.  The following diagram is commutative
$$
\begin{array}{ccc}
S(L_+) & \to & U\A_+ / (U\A_+(\A_+\cap \B)) 
\\ 
\downarrow & \; & \downarrow 
\\
S(L) & \to & 
U\A / (U\A)\B 
\end{array}
$$
The horizontal maps are the linear isomorphisms obtained by
symmetrization. Since the natural map from $S(L_+)$ to $S(L)$ is
injective, so is the right-hand vertical map, and its dual is
surjective.  \hfill \qed\medskip

A {\it Poisson homogeneous space} $(X,P_X)$ for a Poisson-Lie group
$(\Gamma,P_\Gamma)$ is a Poisson formal manifold $(X,P_X)$, equipped
with a transitive action of $\Gamma$, and such that the map
$\Gamma\times X\to X$ is Poisson. Then there is a Lie subgroup $\Gamma'$
of $\Gamma$ such that $X = \Gamma /\Gamma'$. Following \cite{Dr:hom},
$(X,P_X)$ is said to be {\it of group type} if either of the following 
equivalent conditions is satisfied: a) the projection map $\Gamma \to X$
is Poisson, b) the Poisson bivector $P_X$ vanishes at one point of $X$.

\subsection{Poisson homogeneous space structure on $G/H$} 
\label{sect:refinement}

Let $(\G,r,\langle\ ,\ \rangle_{\G},\HH)$ be a quasitriangular Manin
pair.  Let $H$ be the formal subgroup of $G$ with Lie algebra $\HH$,
and let $P_{G/H}$ be the $2$-tensor on $G/H$ equal to the projection of 
$r^{L}$.

\begin{prop} \label{refinement}
Let $(\G,r,\langle\ ,\ \rangle_{\G},\HH)$ be a quasitriangular Manin pair. 
Then $P_{G/H}$ is a Poisson bivector on $G/H$, and $G/H$ is a Poisson
homogeneous space for  $(G,P_G)$.  
\end{prop}

{\em Proof.} The only nonobvious property is the antisymmetry of the
bracket defined by $r^L$. Let $t$ denote the symmetric element of
$\G\otimes\G$ defined by $\langle\ ,\ \rangle_{\G}$.  If $f_1$ and
$f_2$ are right $\HH$-invariant functions on $G$, $\{f_1,f_2\} +
\{f_2,f_1\} = t^{L} ( df_1\otimes df_2) = t^{R}(df_1\otimes df_2)$, by
the invariance of $t$.  Since $\HH$ is Lagrangian, $t$ belongs to $\HH
\otimes \G + \G \otimes \HH$.  Since moreover $f_1$ and $f_2$ are
$\HH$-invariant, $\{f_1,f_2\} + \{f_2,f_1\}$ vanishes.  \hfill
\qed\medskip

\begin{remark} 
In the case where $(\G,r,\langle\ ,\ \rangle_{\G},\HH)$ corresponds to
a Manin quadruple, this statement has been proved by Etingof and Kazhdan,
who also constructed a quantization of this Poisson homogeneous space
(\cite{EK}).
\end{remark}

\begin{remark} \label{rem:gp:type} 
For $g$ in $G$ and $x$ in $\G$, let us denote 
the adjoint action of $g$ on $x$ by $\null^g x$. 
In the case of a quadruple
$(\G,\G_+,\G_-,\HH)$, the Poisson homogeneous space $(G/H,P_{G/H})$ is
of group type if and only if there exists $g$ in $G$ such that $\null^g
\HH$ is {\it graded} for the Manin triple decomposition, i.e., such
that $\null^g \HH = (\null^g \HH \cap \G_+) \oplus (\null^g \HH\cap
\G_-)$.
\end{remark}

\subsection{Poisson homogeneous space structure on $G_+ / (G_+ \cap H)$}

Let $(\G,\G_+,\G_-,\HH)$ be a Manin quadruple.  The inclusion $G_+
\subset G$ induces an inclusion map   $i : G_+ / (G_+ \cap H) \to G/H$. 
On the other hand, by  Proposition \ref{refinement},  
$(\G,\G_+,\G_-,\HH)$ defines a Poisson structure on $G/H$.

\begin{prop} \label{prop:hom:G+}
There exists a unique Poisson structure $P_{G_+ / (G_+ \cap H)}$ on
$G_+ / (G_+ \cap H)$ such that the inclusion $i$ is Poisson.  Then $(
G_+ / (G_+ \cap H), P_{G_+ / (G_+ \cap H)})$ is a Poisson homogeneous
space for $(G_+,P_{G_+})$.
\end{prop}

{\em Proof.}  Let $\langle\ ,\ \rangle_{\G\otimes \G}$ denote 
the bilinear form on $\G\otimes\G$ defined as the tensor square of
$\langle\ ,\ \rangle_{\G}$.  We must prove that for any $g$ in $G_+$,
$r_{\G_+,\G_-}$ belongs to $\G_+ \otimes \G_+ + \G \otimes\null^g\HH +
\null^g\HH \otimes \G$.  The annihilator of this space in
$\G\otimes\G$ with respect to $\langle\ ,\ \rangle_{\G\otimes\G}$ is
$(\G_+\cap \null^g \HH)\otimes \null^g \HH + \null^g\HH \otimes
(\G_+\cap \null^g \HH)$. Let us show that for any $x$ in this
annihilator, $\langle r_{\G_+,\G_-}, x\rangle_{\G\otimes\G}$ is
zero. Assume that $x = v\otimes w$, where $v\in \null^g
\HH, w\in \G_+\cap \null^g \HH$; then
$$
\langle r_{\G_+, \G_-}, x\rangle_{\G\otimes\G} 
= \langle \sum_i e^i \otimes e_i, 
v\otimes  w\rangle_{\G\otimes\G} = \langle v,w\rangle_{\G} = 0, 
$$ 
where the second equality follows from the facts that $\G_+$ is
isotropic and that $(e^i),(e_i)$ are dual bases, and the last equality
follows from the isotropy of $\null^g\HH$. On the other hand, the 
isotropy of $\G_+$ implies that $\langle r_{\G_+,\G_-}, (\G_+\cap \null^g \HH) 
\otimes \null^g \HH \rangle_{\G\otimes\G} = 0$.  Therefore $r_{\G_+,\G_-}$
belongs to $\G_+ \otimes \G_+ + \G \otimes\null^g\HH + \null^g\HH
\otimes \G$; this implies the first part of the Proposition.

Let us prove that $( G_+ / (G_+ \cap H), P_{G_+ / (G_+ \cap H)})$ is a
Poisson homogeneous space for $(G_+,P_{G_+})$. We have a commutative
diagram
$$
\begin{array}{ccc}
G_+ \times (G_+/(G_+\cap H)) & \stackrel{a_{G_+}}{\to} & G_+/(G_+ \cap H)
\\ 
i_G\times i \downarrow & \; & i \downarrow 
\\
G \times (G/H)  & \stackrel{a_G}{\to} & G/H
\end{array}
$$
where $i_G$ is the inclusion map of $G_+$ in $G$ and $a_G$ (resp.,
$a_{G_+}$) is the action map of $G$ on $G/H$ (resp., of $G_+$ on
$G_+ / (G_+ \cap H)$).  The maps $i_G\times i$, $i$ and $a_G$ are
Poisson maps; since $i$ induces an injection of tangent spaces,
it follows that $a_{G_+}$ is a Poisson map.
\hfill\qed\medskip

In \cite{Dr:hom}, Drinfeld defined a Poisson bivector $P'_{G_+ /
(G_+\cap H)}$ on $G_+ / (G_+\cap H)$, which can be described as
follows. When $V$ is a Lagrangian subspace of $\G$, identify $(\G_+\cap
V)^\perp$ with a subspace of $\G_-$, and define $\bar\xi_V : (\G_+\cap
V)^\perp \to \G_+ / (\G_+\cap V)$ to be the linear map which, to any 
element $a_-$ of $(\G_+\cap V)^\perp\subset \G_-$, associates
the class of an element $a_+\in\G_+$ such that $a_+ + a_-$ belongs to
$V$.  Then there is a unique element $\xi_V\in (\G_+ / (\G_+\cap
V))^{\otimes 2}$, such that $(a_- \otimes id)(\xi_V) =
\bar\xi_V(a_-)$, for any $a_-\in (\G_+\cap V)^\perp$.

For any $g$ in $G_+$, identify the tangent space of $G_+ / (G_+\cap
H)$ at $g(G_+\cap H)$ with $\G_+ / (\G_+ \cap
\null^{g}\HH)$ via left-invariant vector fields.  The element
$\xi_{g}$ of $(\G_+ / (\G_+\cap \null^{g}\HH))^{\otimes 2}$
corresponding to the value of $P'_{G_+ / (G_+\cap H)}$ at $g(G_+\cap
H)$ is then $\xi_{\null^{g}\HH}$.

\begin{prop} \label{prop:hom:G+'}
The bivectors $P_{G_+ / (G_+\cap H)}$ and $P'_{G_+ / (G_+\cap H)}$
are equal. 
\end{prop}

{\em Proof.}  We have to show that the injection $i : G_+ / (G_+ \cap H)
\to G/H$ is compatible with the bivectors $P'_{G_+ / (G_+ \cap H)}$ and
$P_{G/H}$. The differential of $i$ at $g(G_+\cap H)$, where $g$ belongs to $G_+$, 
induces the canonical injection $\iota$ from $\G_+ /
(\G_+\cap\null^{g}\HH)$ to $\G / \null^{g}\HH$. 
Let us show that for any $g$ in $G_+$, the injection $\iota\otimes
\iota$ maps $\xi_{g}\in (\G_+ / (\G_+\cap\null^{g}\HH))^{\otimes 2}$
to the class of $r_{\G_+,\G_-}$ in $(\G / \null^{g}\HH)^{\otimes 2}$.
We have to verify the commutativity of the diagram
$$
\begin{array}{ccc}
(\G_+ \cap \null^{g}\HH)^\perp & \stackrel{\iota^*}{\leftarrow} & 
(\null^{g}\HH)^\perp
\\ 
\bar\xi_{g} \downarrow & \; & \bar r_{g} \downarrow 
\\
\G_+ / (\G_+\cap\null^{g}\HH) & \stackrel{\iota}{\rightarrow} & \G 
/ \null^{g}\HH
\end{array}
$$
where the horizontal maps are the natural injection and restriction
maps, and $\bar r_{g}$ is defined by $\bar r_{g}(a) = $ the class of
$\langle r_{\G_+,\G_-}, a \otimes id\rangle_{\G\otimes\G}$ mod
$\null^{g}\HH$, for $a\in (\null^{g}\HH)^\perp$. Let $a$ belong to
$(\null^{g}\HH)^\perp$.  By the maximal isotropy of
$\null^{g}\HH$, the element $a$ can be identified with an element of
$\null^{g}\HH$.  Let us write $a = a_+ + a_-$, with
$a_\pm\in\G_\pm$.  Then $\iota^*(a) = a_-$, $\bar\xi_{g}(a_-) = a_+
+ (\G_+ \cap\null^{g}\HH)$ by the definition of $\bar\xi_{g}$ and
because $a\in\null^{g}\HH$, and $\iota(a_+ + \G_+
\cap\null^{g}\HH) = a_+ + \null^{g}\HH$. On the other hand, $\bar
r_{g}(a) = a_+ + \null^{g}\HH$, so the diagram commutes.  \hfill
\qed\medskip

Moreover, it is a result of Drinfeld (\cite{Dr:hom}, Remark 2)  that the
formal Poisson homogeneous spaces for $(G_+,P_{G_+})$ are all of the
type described in Proposition \ref{prop:hom:G+}.

\begin{remark} The $(G_+,P_{G_+})$-Poisson homogeneous space  
$( G_+ / (G_+ \cap H), P_{G_+ / (G_+ \cap H)})$ is of group type if
and only if for some $g\in G_+$, $\null^{g}\HH$ is graded for the
Manin triple decomposition, see Remark \ref{rem:gp:type}.
\end{remark}

\begin{remark} {\it Conjugates of Manin quadruples.}
Let us denote the conjugate of an element $x$ in $G$ by  an element
$h$ in $G$ by $\null^g x = g x g^{-1}$. For  $(\G,\G_+,\G_-,\HH)$ a Manin
quadruple, and $g$ in $G$,  $(\G,\G_+,\G_-,\null^g \HH)$ is also a Manin
quadruple. If $g$ belongs to $H$,  this is the same quadruple; and if
$g$ belongs to $G_+$, the Poisson structure induced by 
$(\G,\G_+,\G_-,\null^g \HH)$ on $G_+ / (G_+ \cap \null^g H)$  is
isomorphic to that induced by  $(\G,\G_+,\G_-,\HH)$ on $G_+ / (G_+ \cap
H)$, via conjugation by $g$.  It follows that Poisson homogeneous space
structures induced by  $(\G,\G_+,\G_-,\HH)$ are indexed by elements of
the double quotient $G_+ \setminus  G / H$.  
\end{remark}

\begin{remark}
The first examples of Poisson homogeneous spaces which are not of
group type were studied in \cite{Dazord} and \cite{Lu}, under the name
of affine Poisson structures. In these cases, the stabilizer of a point 
is trivial.  
\end{remark}

\section{Quantization of $G/H$}

In this section, we introduce axioms for the quantization of a
quasitriangular Manin pair. We then show that any solution to this
quantization problem leads to a quantization of  the Poisson homogeneous
space $G/H$ constructed in Section \ref{sect:refinement}.

\subsection{Definition of quantization of Poisson homogeneous spaces}

Let $(\Gamma,P_\Gamma)$ be a Poisson-Lie group and let $(X,P_X)$ be a
formal Poisson homogeneous space over $(\Gamma,P_\Gamma)$. Let
$(\cA,\Delta_\cA)$ be a quantization of the enveloping algebra of the
Lie algebra of $\Gamma$, and let $\cA^{opp}$ denote the opposite
algebra to $\cA$.

\begin{defin}
A {\em quantization of the Poisson homogeneous space}
$(X,P_X)$ is a $\CC[[\hbar]]$-algebra $\cX$, such that

1) $\cX$ is a quantization of the Poisson algebra $(\cO_X,P_X)$
of formal functions on $X$, and 

2) $\cX$ is equipped with an algebra-module structure over
$(\cA^{opp},\Delta)$, whose reduction mod $\hbar$ coincides with the
algebra-module structure of $\cO_X$ over $((U\G)^{opp},\Delta_0)$, 
where $\Delta_0$ is the coproduct of $U\G$. 
\end{defin}

Condition 1 means that there is an isomorphism
of $\CC[[\hbar]]$-modules from $\cX$ to $\cO_X[[\hbar]]$, inducing an 
algebra isomorphism between $\cX / \hbar \cX$ and $\cO_X$, and inducing 
on $\cO_X$ the Poisson structure defined by $P_X$.    

The first part of condition 2 means that $\cX$ has a module structure
over $\cA^{opp}$, such that for $x,y\in \cX, a\in \cA$, and $\Delta(a)
= \sum a^{(1)} \otimes a^{(2)}$, $a(xy) = \sum a^{(1)}(x) a^{(2)}(y)$.

\subsubsection*{Conventions}  We will say that a $\CC[[\hbar]]$-module
$V$ is {\it topologically free} if it is  isomorphic to $W[[\hbar]]$, where
$W$ is a complex vector space. We denote  the canonical projection of
$V$ onto $V / \hbar V$ by  $v\mapsto v$ mod $\hbar$.    In what follows,
all tensor products of $\CC[[\hbar]]$-modules are $\hbar$-adically
completed. When $E$ is a $\CC[[\hbar]]$-module, we denote by $E^*$ its dual
$\Hom_{\CC[[\hbar]]}(E,\CC[[\hbar]])$.  For a subset $\cS$ of an algebra
$\cA$, we denote by $\cS^\times$ the group of invertible elements of
$\cS$. When $\cA,\cB$ are two Hopf or quasi-Hopf algebras with unit elements
$1_\cA,1_\cB$ and counit maps $\eps_\cA,\eps_\cB$, and $\cS$ is a
subset of $\cA \otimes\cB$, we denote by $\cS^\times_0$ the subgroup of
$\cS^\times$ with elements $x$ such that $(\eps_\cA\otimes id)(x) =
1_\cB$ and $(id \otimes \eps_\cB)(x) = 1_\cA$. We also denote
$\Ker\eps_\cA$ by $\cA_0$.

\subsection{Quantization of quasitriangular Manin pairs}

\begin{defin} \label{defin:quant:1}
A {\em quantization of a quasitriangular Manin pair} $(\G,r,\langle\ ,\
\rangle_\G,\HH)$ is the data of 

\noindent 1) a quasitriangular Hopf algebra $(A,\Delta,\cR)$  quantizing
$(\G,r,\langle\ ,\ \rangle_\G)$,  

\noindent 2) a subalgebra $B\subset A$ and an element $F$ in $(A\otimes
A)_0^\times$, such that  

a) $B\subset A$ is a flat deformation of the inclusion $U\HH\subset U\G$,  

b) $F\Delta(B)F^{-1} \subset B\otimes B$, and  
$F^{(12)} (\Delta\otimes id)(F)  \left( F^{(23)} (id
\otimes \Delta)(F) \right)^{-1}\in B^{\otimes 3}$, 

c) there exists a Lagrangian complement $L$ of $\HH$ in $\G$, such that
$$
\left( {1\over\hbar}(F - F^{(21)}) \ mod\ \hbar \right)  
= r - r_{\HH,L},
$$ 
where $r_{\HH,L} = \sum_i \eps^i\otimes \eps_i$, and 
$(\eps^i),(\eps_i)$ are dual bases of $\HH$ and $L$.  
\end{defin}

Observe that condition 1 implies that $( {1\over\hbar}(\cR - 1)$ mod
$\hbar) = r$.  

In condition 2, $(A\otimes A)^\times_0 = 1 + \hbar (A_0\otimes A_0)$.

Condition 2a means that the $\CC[[\hbar]]$-module isomorphism
between $A$ and $U\G[[\hbar]]$ arising from condition 1 can be chosen
in such a way that it induces an isomorphism between $B$ and
$U\HH[[\hbar]]$.

Let $(A,\Delta_B,\Phi_B)$ be the quasi-Hopf algebra obtained by twisting
the Hopf algebra $(A,\Delta)$ by $F$. Then, by definition, $\Delta_B(x)
= F \Delta(x) F^{-1}$, for any $x\in A$, and $\Phi_B = F^{(12)}
(\Delta\otimes id)(F) \left( F^{(23)} (id \otimes \Delta)(F)
\right)^{-1}$.  Condition 2b expresses the fact that $B\subset A$ is a
sub-quasi-Hopf algebra of $A$.

Condition 2c expresses the fact that the classical limit of
$(B,\Delta_B,\Phi_B)$ is $\HH$ equipped with the Lie quasi-bialgebra
structure associated with $L$.

\subsection{Quantization of $(G/H,P_{G/H})$}

Let us denote the counit map of $A$ by $\eps$.  

\begin{thm} \label{thm:main}
Assume that $((A,\Delta,\cR),B,F)$ is a quantization of the 
quasitriangular Manin pair $(\G,r,\langle\ ,\ \rangle_\G,\HH)$.  
Let $(A^*)^B$ be the subspace of $A^*$ consisting of the forms $\ell$
on $A$ such that $\ell(ab) = \ell(a)\eps(b)$ for any $a\in A$ and
$b\in B$.  

a) For $\ell,\ell'$ in $(A^*)^B$, define $\ell * \ell'$
to be the element of $A^*$ such that 
$$
(\ell * \ell')(a) = (\ell \otimes  \ell')(\Delta(a) F^{-1}) , 
$$
for any $a$ in $A$. Then $*$ defines an associative algebra structure on
 $(A^*)^B$. 

b) For $a\in A$ and $\ell \in (A^*)^B$, define~~ $a\ell$ ~~to be the form
on $A$ such that $(a\ell)(a') = \ell(aa')$, for any $a'\in A$. This
map defines on $((A^*)^B,*)$ a structure of an algebra-module over the
Hopf algebra $(A^{opp},\Delta)$.

c) The algebra $((A^*)^B,*)$ is a quantization of the Poisson algebra
$(\cO_{G/H} , P_{G/H})$.

With its algebra-module structure over $(A^{opp},\Delta)$,
$((A^*)^B,*)$ is a quantization of the Poisson homogeneous
space $(G/H,P_{G/H})$.
\end{thm}

{\em Proof.} Let $\ell$ and $\ell'$ belong to $(A^*)^B$. Then 
for any $a\in A,b\in B$,  
\begin{align*}
& (\ell*\ell')(ab) = (\ell \otimes\ell')(\Delta(a)\Delta(b) F^{-1})  
=  (\ell \otimes\ell')(\Delta(a) F^{-1} \Delta_B (b) )  
\\ & = \eps(b) 
(\ell \otimes\ell')(\Delta(a) F^{-1} )
=  \eps(b)  (\ell*\ell')(a) , 
\end{align*}
where the third equality follows from the fact that $\Delta_B(B)
\subset B\otimes B$ and $(\eps \otimes\eps) \circ\Delta_B = \eps$.  It
follows that $\ell*\ell'$ belongs to $(A^*)^B$.

Let $\ell,\ell'$ and $\ell''$ belong to $(A^*)^B$. Then for any $a$ in
$A$, 
\begin{equation} \label{maoz}
((\ell * \ell') *\ell'')(a) = (\ell \otimes\ell' \otimes \ell'') \left( 
(\Delta\otimes id)\circ\Delta(a) (\Delta\otimes id)(F^{-1}) (F^{(12)})^{-1} \right) , 
\end{equation}
and 
$$
(\ell * (\ell' *\ell''))(a) = (\ell \otimes\ell' \otimes \ell'') \left(  (id
\otimes \Delta)\circ\Delta(a) (id \otimes \Delta)(F^{-1})
(F^{(23)})^{-1} \right) . 
$$
By the coassociativity of $\Delta$ and the definition of $\Phi_B$, this
expression is equal to $(\ell \otimes\ell' \otimes \ell'') \left( 
(\Delta \otimes id)\circ\Delta(a) (\Delta\otimes id)(F^{-1})
(F^{(12)})^{-1}\Phi_B^{-1} \right)$, and since $\Phi_B$ belongs to
$B^{\otimes 3}$ and $\eps^{\otimes 3}(\Phi_B) = 1$, this is equal to the
right-hand side of (\ref{maoz}). It follows that $(\ell * \ell') *
\ell'' = \ell * (\ell' * \ell'')$, so $*$ is associative. Moreover,
$\eps$ belongs to $(A^*)^B$ and is the unit element of $((A^*)^B,*)$.
This proves part a of the theorem. 

Part b follows from the definitions. 

Let us prove that $(A^*)^B$ is a flat deformation of $(U\G / (U\G)\HH
)^*$. Let us consider a Lagrangian complement $L$ of $\HH$ in $\G$. 
Let
$\Sym$ denote the symmetrisation map from the symmetric algebra
$S(\G)$ of $\G$ to $U\G$, i.e., the unique linear map such that
$\Sym(x^l) = x^l$ for any $x\in\G$ and $l\geq 0$, and let us define
$\wt L$ to be $\Sym(S(L))$. Thus $\wt L$ is a linear subspace of
$U\G$, and inclusion of  $\wt L\otimes U\HH$ in $U\G\otimes U\G$
followed by multiplication induces a linear
isomorphism from $\wt L\otimes U\HH$ to $U\G$. It follows that the
restriction to $\wt L$ of the projection $U\G \to U\G /(U\G)\HH$ is an
isomorphism, which defines a linear isomorphism between $\cO_{G/H} = (U\G /
(U\G)\HH )^*$ and $\wt L^*$.

Let us fix an isomorphism of $\CC[[\hbar]]$-modules from $A$ to
$U\G[[\hbar]]$, inducing an isomorphism between $B$ and $U\HH[[\hbar]]$ 
 and let us define $C$ to be the preimage of $\wt L[[\hbar]]$. Thus  
$C$ is
isomorphic to $\wt L[[\hbar]]$, and inclusion followed by multiplication
induces a linear isomorphism between $C\otimes B$ and $A$, as does any
morphism between two topologically free  $\CC[[\hbar]]$-modules $E$ and
$F$, which induces an isomorphism between $E / \hbar E$ and $F / \hbar
F$. Therefore, restriction of linear forms to $C$ induces an
isomorphism between $(A^*)^B$ and $C^*$. It follows that $(A^*)^B$ is
isomorphic to $(\wt L[[\hbar]])^*$, which is in turn isomorphic  to 
$\wt L^*[[\hbar]]$. Since $\wt L^*$  is isomorphic to $\cO_{G/H}$, 
$(A^*)^B$ is isomorphic to $\cO_{G/H}[[\hbar]]$. 

Let us fix $\ell,\ell'$ in $(A^*)^B$, and let us compute $({1\over \hbar}
(\ell * \ell' - \ell' * \ell)$ mod $\hbar)$.  Let us set 
$f = ({1\over\hbar}(F-1)$ mod $\hbar)$.  
For $a$ in $A$,  
\begin{align*}
& {1\over \hbar} (\ell * \ell' - \ell' * \ell)(a) = 
(\ell \otimes\ell') \left( {1\over \hbar} ( \Delta(a) F^{-1} - 
\Delta'(a) (F^{(21)})^{-1} ) \right)
\\ & 
= (\ell \otimes\ell') \left( {1\over \hbar} (\Delta(a) F^{-1} - 
\cR \Delta(a) \cR^{-1} (F^{(21)})^{-1} ) \right)
\\ & = 
(\ell \otimes \ell') \left( - r \Delta_0(a_0) + \Delta_0(a_0)(r 
+ f^{(21)} - f) \right) + o(\hbar),  
\end{align*}
where $\Delta_0$ is the coproduct of $U\G$ and $a_0$ is the image of
$a$ in $A / \hbar A = U\G$. 
Since, by condition 2c of Definition \ref{defin:quant:1},
$r - f + f^{(21)}$ is equal to $r_{\HH,L}$, it belongs to $\HH\otimes\G$ 
and since $\ell$ and $\ell'$ are right $B$-invariant,   
$({1\over \hbar} (\ell * \ell' - \ell' * \ell)(a)$ mod $\hbar) =  
(\ell \otimes \ell')( - r \Delta_0(a_0))$, 
which is the Poisson bracket defined by $r^{L}$ on $G/H$.
This ends the proof of part c of the theorem.   

Since $(F$ mod $\hbar) = 1$, the reduction modulo $\hbar$ of the
algebra-module  structure of $((A^*)^B,*)$ over $(A^{opp},\Delta)$ is
that of $\cO_X$ over $((U\G)^{opp},\Delta_0)$. This ends the proof of
the theorem.   \hfill \qed\medskip 

\begin{remark}  If $((A,\Delta,\cR),B,F)$ is a quantization of a
quasitriangular Manin pair $(\G,r,\langle\ ,\ \rangle_\G,\HH)$, and if
$F_0$ is an element of  $(B\otimes B)^\times_0$, then 
$((A,\Delta,\cR),F_0F)$ is a quantization  of the same Manin pair. We
observe that the product $*$ on $(A^*)^B$ is independent of such a
modification of $F$. 
\end{remark}

\begin{remark} \label{rem:graded} 
In the case of a Manin quadruple $(\G,\G_+,\G_-,\HH)$ where $\HH$ is
graded for the Manin triple decomposition (see Remark \ref{rem:gp:type}),  
$r_{\G_+,\G_-} - r_{\HH,L} = f - f^{(21)}$ belongs to $\HH\otimes\G 
+ \G\otimes\HH$. The corresponding quantum condition is that 
\begin{equation} \label{condition:F}
F \in 1 + \hbar (B_0\otimes A_0 + A_0 \otimes B_0).
\end{equation}
This is the case is \cite{ER}, where $F$
belongs to $(B\otimes A)^\times_0 = 1 + \hbar (B_0 \otimes A_0)$.   

When condition (\ref{condition:F}) is fulfilled,  the product $*$ is
the restriction to $(A^*)^B$ of the usual product on  $A^*$, defined as
the dual map to $\Delta$.     
\end{remark}

\subsection{Relations in $(A^*)^B$}

It is well-known that the matrix coefficients of the representations of a
quasitriangular Hopf algebra can be organized in $L$-operators, 
satisfying the so-called $RLL$ relations.  We recall this construction 
and introduce analogues of these matrix coefficients and of the $RLL$
relations  for the algebra $(A^*)^B$. 

Recall that there is an algebra structure on $A^*$, where the  product
is $(\ell,\ell') \mapsto \ell\ell'$, such that for any $a\in A$,
$(\ell\ell')(a) = (\ell\otimes\ell')(\Delta(a))$. 

Let $\Rep(A)$ be the category of modules over $A$, which are free and
finite-dimensional over $\CC[[\hbar]]$.  There is a unique map
$$
\oplus_{V\in \Rep(A)} (V^*\otimes V) \to A^*, \quad 
\kappa \mapsto \ell_\kappa,  
$$
such that for any object $(V,\pi_V)$ in $\Rep(A)$, and any $\xi\in
V^*$ and $v\in V$, $\ell_{\xi\otimes v}(a) = \xi(\pi_V(a)v)$, for any
$a\in A$. Define $\Coeff(A)$ to be the image of this map. Then
$\Coeff(A)$ is a subalgebra of $A^*$. Moreover, there is a Hopf
algebra structure on $\Coeff(A)$, with coproduct $\Delta_{\Coeff(A)}$
and counit $\eps_{\Coeff(A)}$, uniquely determined by the rules
$$
\Delta_{\Coeff(A)}(\ell_{\xi\otimes v}) = \sum_{i} \ell_{\xi\otimes
v_i} \otimes \ell_{\xi^i \otimes v}, \quad
\eps_{\Coeff(A)}(\ell_{\xi\otimes v}) = \xi(v),
$$
where $(v_i)$ and $(\xi^i)$ are dual bases of $V$ and $V^*$. The
duality pairing between $A$ and $A^*$ then induces a Hopf algebra
pairing between $(A,\Delta)$ and $(\Coeff(A),\Delta_{\Coeff(A)})$ (see
\cite{Andr}).

For $V$ an object of $\Rep(A)$, define $L_V$ to be  the element of
$\End(V)\otimes \Coeff(A)$ equal to  $\sum_{i} \kappa^i \otimes
\ell_{\kappa_i}$, where  $(\kappa^i)$ and $(\kappa_i)$ are dual bases of
$\End(V)$  and  $V^*\otimes V$.  It follows from $\cR \Delta = \Delta' \cR$ 
that the relation 
$$
R^{(12)}_{V,W} L_V^{(1a)} L_W^{(2a)} = 
L_W^{(2a)} L_V^{(1a)} R_{V,W}^{(12)}
$$
is satisfied in $\End(V)\otimes \End(W) \otimes \Coeff(A)$, where the
superscripts $1,2$ and $a$ refer to the successive factors of
the tensor product. Moreover, 
$$
(id_{V}\otimes \Delta_{\Coeff(A)}) (L_V) = 
L_V^{(1a)} L_V^{(1a')}
$$  
holds in 
$\End(V)\otimes \Coeff(A)^{\otimes 2}$, where the superscripts 
$1,a$ and $a'$ refer to the successive factors of this tensor product. 

For $(V,\pi_V)$ an $A$-module,  we set $V^B = \{v\in V | \forall b\in B,
\pi_V(b)(v) =  \eps(b)v\}$. There is a unique map 
$$
\oplus_{V\in \Rep(A)} (V^*\otimes V^B) \to (A^*)^B, 
\quad \kappa \mapsto \wt\ell_\kappa
$$
such that for $\xi\in V^*$ and $v\in V^B$, $\wt\ell_{\xi\otimes v}(a)
= \xi(\pi_V(a)v)$, for any $a\in A$.  Define $\Coeff(A,B)$ to be the image 
of this map. 

For $V$ an object of $\Rep(A)$, $V^B$ is a free, finite dimensional
$\CC[[\hbar]]$-module. It follows that the dual of $V^*\otimes V^B$ is
$(V^B)^* \otimes V$, which may be identified with $\Hom_{\CC[[\hbar]]}
(V^B,V)$. Define $\wt L_V$ to be the element of $\Hom_{\CC[[\hbar]]}
(V^B,V) \otimes \Coeff(A,B)$ equal to $\sum_i \kappa^i \otimes
\wt\ell_{\kappa_i}$, where $(\kappa^i)$ and $(\kappa_i)$ are dual
bases of $\Hom_{\CC[[\hbar]]} (V^B,V)$ and $V^*\otimes V^B$.

When $(V,\pi_V)$ and $(W,\pi_W)$ are objects of $\Rep(A)$, let
$R_{V,W}$ be the element of $\End_{\CC[[\hbar]]}(V\otimes W)$ equal to
$(\pi_V\otimes \pi_W)(\cR)$, where $\cR$ is the $R$-matrix of $A$.
Recall that the twist of $\cR$ by $F$ is $\cR_B = F^{(21)} \cR
F^{-1}$, and set $R_{B; V,W} = (\pi_V\otimes \pi_W)(\cR_B)$.

\begin{prop} 
$\Coeff(A,B)$ is a subalgebra of $(A^*)^B$.  For any objects $V$ and
$W$ in $\Rep(A)$, the relation 
$$ 
R^{(12)}_{V,W} \wt L^{(1a)}_V \wt L^{(2a)}_W = 
\wt L^{(2a)}_W \wt L^{(1a)}_V (R^{(12)}_{B;V,W})_{| Z}
$$
is satisfied in $\Hom_{\CC[[\hbar]]} (Z,V\otimes W) \otimes
\Coeff(A,B)$, where $Z$ is the intersection $(V^B\otimes W^B)\cap
R_{B;V,W}^{-1}(V^B\otimes W^B)$.  In this equality, the left-hand side
is an element of $\Hom_{\CC[[\hbar]]}(V^B\otimes W^B, V\otimes W)
\otimes \Coeff(A,B)$, viewed as an element of
$\Hom_{\CC[[\hbar]]}(Z, V\otimes W) \otimes \Coeff(A,B)$ by
restriction.
\end{prop}

Recall that an algebra-comodule $\cX$ over a Hopf algebra
$(\cA,\Delta_\cA)$ is the data of an algebra structure over $\cX$ and
a left comodule structure of $\cX$ over $(\cA,\Delta_\cA)$,
$\Delta_{\cX,\cA} : \cX \to \cA \otimes \cX$, which is also a morphism
of algebras.

\begin{prop}
There is a unique algebra-comodule structure on $\Coeff(A,B)$
over $(\Coeff(A),\Delta_{\Coeff(A)})$, compatible with the algebra-module 
structure of $(A^*)^B$ over $(A^{opp},\Delta)$.  The relation  
$$
(id_V \otimes \Delta_{ \Coeff(A,B), \Coeff(A) })(\wt L_V) 
= L_V^{(1a)} \wt L_V^{(1a')} 
$$ 
is satisfied in $\End(V) \otimes \Coeff(A) \otimes \Coeff(A,B)$, where
the superscripts  $1,a$ and $a'$ refer to the successive factors of
this tensor product.   
\end{prop}

\section{Quantization of $G_+ / (G_+\cap H)$}

In this section, we state axioms for the quantization of a Manin 
quadruple, and show that any such quantization gives rise to a
quantization of the Poisson homogeneous space $G_+ / (G_+ \cap H)$
constructed in \cite{Dr:hom} (see Propositions \ref{prop:hom:G+} 
and \ref{prop:hom:G+'}).  

\subsection{Quantization of Manin quadruples}

Let us fix a Manin quadruple $(\G,\G_+,\G_-,\HH)$. Recall that  
$(\G,\G_+,\G_-,\HH)$ gives rise to a quasitriangular Manin pair  
$(\G,r,\langle\ ,\ \rangle_\G,\HH)$, if we set 
$r = r_{\G_+,\G_-} = \sum_i e^i \otimes e_i$, where
$(e^i)$ and $(e_i)$ are dual bases of $\G_+$ and $\G_-$.    

\begin{defin} 
A quantization of a Manin quadruple $(\G,\G_+,\G_-,\HH)$ is the data of 

\noindent 1) a quantization $((A,\Delta,\cR),B,F)$ of 
the quasitriangular Manin pair 
$(\G,r_{\G_+,\G_-},\langle\ ,\  \rangle_\G,\HH)$,  

\noindent 2) a Hopf subalgebra  $A_+$ of $(A,\Delta)$
such that  

a) $A_+ \subset A$ is a flat deformation of  $U\G_+\subset U\G$,

b) $B\cap A_+\subset A_+$ is a flat deformation of the inclusion 
$U(\HH\cap\G_+)\subset U\G_+$, 

c) $F$ satisfies  
\begin{equation} \label{hyp} 
F^{-1} \in \left( (AB_0 + A_+)\otimes A + A\otimes AB_0 \right) \cap  
\left( AB_0 \otimes A + A\otimes (AB_0 + A_+) \right) .  
\end{equation} 
\end{defin}

It follows from
condition 2c of Definition \ref{defin:quant:1} and the beginning of the
proof of  Proposition \ref{prop:hom:G+} that $({1\over\hbar}(F -
F^{(21)})$  mod $\hbar)$ belongs to $\G_+ \otimes\G_+ + \HH\otimes\G + \G
\otimes\HH = \left( (\HH + \G_+) \otimes\G + \G\otimes \HH \right)  \cap
 \left( \HH \otimes\G + \G\otimes  (\HH + \G_+) \right)$. Therefore condition 
(\ref{hyp}) is natural. It is equivalent to the condition that $F$ belong 
to the product of subgroups of $(A\otimes A)^\times_0$ 
$$
\left( 1 + \hbar (AB_0 \otimes A_0 + A_0 \otimes AB_0) \right) 
\left( 1 + \hbar (A_+)_0 \otimes (A_+)_0 \right) . 
$$

\subsubsection*{Example}

Recall that (\ref{infinite:quadruple}) is a graded Manin quadruple. 
In  \cite{ER}, a quantization  of this quadruple was constructed 
for the case where $\bar\G = \SL_2$. 

\begin{remark} 
In the case where $\G_+ \cap \HH = 0$, which corresponds  to a  
homogeneous space over $G_+$ with trivial stabilizer,  $F$ automatically
satisfies condition  (\ref{hyp}). Indeed, in that case, multiplication
induces an  isomorphism $A_+ \otimes B \to A$, therefore $A = AB_0 +
A_+$.  
\end{remark}

\subsection{Quantization of $(G_+ / (G_+ \cap H) ,  P_{G_+ / (G_+ \cap H) })$}

Let $I_0$ be the subspace of $\cO_{G/H}$ equal to  $\cO_{G/H} \cap
(U\G_+)^\perp$. It follows from Lemma \ref{lemma:surj} that  $I_0$ is an
ideal of $\cO_{G/H}$, and that the algebra  $\cO_{G_+ / (G_+\cap H)}$
can be identified as a Poisson algebra with the quotient $\cO_{G/H} / I_0$.

Let $I$ be the subspace of $(A^*)^B$ defined as the set of all 
linear forms $\ell$ on $A$ such that $\ell(a_+) = 0$ for any $a_+\in A_+$. 
Therefore 
$$
I = (A^*)^B \cap A_+^\perp . 
$$

\begin{thm} \label{thm:bis} 
Assume that $((A,\Delta,\cR),A_+,B,F)$ is a quantization of the Manin
quadruple $(\G,\G_+,\G_-,\HH)$.  Then $I$ is a two-sided ideal in
$(A^*)^B$, the algebra $(A^*)^B / I$ is a flat deformation of
$\cO_{G_+ / (G_+ \cap H)}$, and is a quantization of the formal
Poisson space $(G_+ / (G_+ \cap H), P_{G_+ / (G_+ \cap H)})$.

Moreover, $I$ is preserved by  the action of $A_+^{opp}$, and 
$(A^*)^B / I$
is an algebra-module over $(A_+^{opp},\Delta)$. With this
algebra-module structure,  $(A^*)^B / I$ is a quantization of the
$(G_+,P_{G_+})$-Poisson homogeneous space $(G_+ / (G_+ \cap H),  P_{G_+ /
(G_+ \cap H)})$.  
\end{thm}

{\em Proof.} Let us fix $\ell$ in $I$ and $\ell'$ in $(A^*)^B$. 
For any $a_+$ in $A_+$, we have 
$$
(\ell * \ell')(a_+) = (\ell\otimes \ell')( \Delta(a_+) F^{-1}) = 0
$$
because $\Delta(A_+) \subset A_+\otimes A_+$ and by assumption 
(\ref{hyp}) on $F$. In the same way, $(\ell' * \ell)(a_+) = 0$. Therefore
$I$ is a two-sided ideal in $(A^*)^B$. 

Let us prove that $(A^*)^B / I$ is a flat deformation of $\cO_{G_+ /
(G_+ \cap H)}$.  To this end, we will identify the
$\CC[[\hbar]]$-modules $(A^*)^B / I$ with $(A_+^*)^{A_+\cap B}$, to
which we apply the result of Theorem \ref{thm:main}.

Recall that $B_0$ and $(A_+\cap B)_0$ denote the augmentation ideals of
$B$ and $A_+\cap B$. Thus $(A^*)^B$ is equal to  $(A / A B_0)^*$, where
$AB_0$ is the image of the product map  $A\otimes B_0\to A$.  In the
same way, $(A_+^*)^{A_+\cap B}$ is equal to  $(A_+ / A_+(A_+\cap
B)_0)^*$.  Let us show that $(A^*)^B/I$ is equal to  $(A_+ / A_+(A_+\cap
B)_0)^*$.  Restriction of a linear form to $A_+$ induces a linear map
$\rho :  (A / A B_0)^* \to (A_+ / A_+(A_+\cap B)_0)^*$. Moreover, the
kernel of $\rho$ is $I$, therefore $\rho$ induces an injective map 
$$
\wt\rho: (A / A B_0)^* / I \to (A_+ /
A_+(A_+\cap B)_0)^*. 
$$

Let us now show that $\wt\rho$ is surjective. For this, it is  enough to
show that the restriction map  $\rho : (A / A B_0)^* \to (A_+ /
A_+(A_+\cap B)_0)^*$ is surjective. $(A / A B_0)^*$ and $(A_+ /
A_+(A_+\cap B)_0)^*$ are topologically free $\CC[[\hbar]]$-modules,  and
the map from  $(A / A B_0)^*   / \hbar   (A / A B_0)^*$
to $(A_+ / A_+(A_+\cap B)_0)^*  / \hbar (A_+ / A_+(A_+\cap
B)_0)^*  $ coincides with the canonical map from $\cO_{G/H}$ to 
$\cO_{G_+ / (G_+ \cap H)}$ which, by Lemma \ref{lemma:surj},  is
surjective. Therefore $\rho$ is surjective, and so is $\wt\rho$.  It
follows that $(A_+^*)^{A_+\cap B}$ is a flat deformation of  $\cO_{G_+ /
(G_+\cap H)}$. 

There is a commutative diagram of algebras 
$$
\begin{array}{ccc}
(A^*)^B
 & \to & (A^*)^B / I 
\\ 
\downarrow & \; & \downarrow 
\\
\cO_{G/H}  & \to & \cO_{G_+ / (G_+\cap H)}
\end{array}
$$
where the vertical maps are projections $X\to X / \hbar X$. 
Since the projection $\cO_{G/H} \to \cO_{G_+/(G_+\cap H)}$ is a 
morphism of Poisson algebras, where   $\cO_{G/H}$ and $\cO_{G_+/(G_+\cap
H)}$ are equipped with  $P_{G/H}$ and $P_{G_+ / (G_+ \cap H)}$, the
classical limit of   $(A^*)^B / I$ is $(\cO_{G_+/(G_+\cap H)},$ $P_{G_+ /
(G_+ \cap H)})$. 

Finally, the algebra-module structure of $(A^*)^B$ over
$(A^{opp},\Delta)$ induces by restriction an algebra-module structure
on $(A^*)^B$ over $(A_+^{opp},\Delta)$, and since $I$ is preserved by
the action of $A_+^{opp}$, $(A^*)^B / I$ is also an algebra-module
over $(A_+^{opp},\Delta)$.  \hfill \qed\medskip

\begin{remark} \label{rem:conj}
For $u\in A^\times$, set $\null^u B = u B u^{-1}$ and 
$\null^u F = (u\otimes u)
\Delta_B(u)^{-1} F$. Let $((A^*)^{\null^u B},*_u)$ be the 
algebra-module over $(A^{opp},\Delta)$   corresponding to 
$(A,\null^u B,\null^u F)$. There is an algebra-module isomorphism 
$i_u : ((A^*)^{B},*) \to ((A^*)^{\null^u B},*_u)$, given by $(i_u\ell)(a) = \ell(au)$, 
for any $a\in A$. 

If $u$ does not lie in $A_+$, there is no reason for  $A_+ \cap \null^u
B$ to be a flat deformation of its classical limit, nor for  $\null^u F$
to satisfy  (\ref{hyp}). But, if the conditions of Theorem \ref{thm:bis}
are still valid for  $(A,\null^u B,\null^u F)$, the resulting algebra-module
over $(A_+^{opp},\Delta_+)$ can be different from the one arising from
$(A,B,F)$. We will see an example of this  situation in section
\ref{sect:spheres}.  \end{remark}

\begin{remark}
In their study of preferred deformations, Bonneau {\it et al.}\ 
studied the case of quotients of compact, connected Lie groups
(\cite{Flato}).  

In \cite{DGS}, Donin, Gurevich and Shnider used quasi-Hopf algebra
techniques to construct quantizations of some homogeneous
spaces. More precisely, they classified the Poisson homogeneous structures
on the semisimple orbits of a simple Lie group with Lie algebra
$\G_0$, and constructed their quantizations using Drinfeld's series
$F_\hbar$ relating the Hopf algebra $U_\hbar \G_0$ to a quasi-Hopf
algebra structure on $U\G_0[[\hbar]]$ involving the Knizhnik-Zamolodchikov 
associator. 

In \cite{Par}, Parmentier also used twists to propose a quantization scheme 
of Poisson structures on Lie groups, generalizing the affine Poisson 
structures. 
\end{remark}

\section{Examples} \label{sect:examples}

In view of Remark \ref{rem:graded}, we can only find nontrivial
applications of the above results in the case of a nongraded $\HH$. In
this section, we shall construct quantizations of some nongraded Manin
quadruples.

\subsection{Finite dimensional examples} \label{sect:spheres}

Let us set $\bar\G = \SL_2(\CC)$; let $\bar\G = \bar\N_+ \oplus\bar\HH
\oplus\bar\N_-$ be the Cartan decomposition of $\bar\G$, and let
$(\bar e_+,\bar h,\bar e_-)$ be the Chevalley basis of $\bar\G$, so
$\bar\N_\pm = \CC e_\pm$ and $\bar\HH = \CC \bar h$. Let $\langle\ ,\
\rangle_{\bar\G}$ be the invariant symmetric bilinear form on $\bar\G$
such that $\langle \bar h ,\bar h \rangle_{\bar\G} = 1$. Set $\G =
\bar\G \times \bar\G$, $\langle (x,y), (x',y') \rangle_{\G} = \langle
x,x' \rangle_{\bar\G} - \langle y,y' \rangle_{\bar\G}$. Set $\G_+ =
\{(x,x), x\in\bar\G\}$ and $\G_- = \{(\eta + \xi_+, - \eta + \xi_- ),
\xi_\pm \in \bar\N_\pm, \eta\in\bar\HH\}$. Then $(\G,\G_+,\G_-)$ is a
Manin triple. Let $(e_+,h,e_-,e_+^*,h^*,e_-^*)$ be the basis of $\G$,
such that $x = (\bar x,\bar x)$ for $x\in\{e_+,h,e_-\}$, $e_+^* =
(\bar e_+,0)$, $h^* = (\bar h,-\bar h)$ and $e_-^* = (0,\bar e_-)$.

A quantization of the Manin triple $(\G,\G_+,\G_-)$ is the algebra $A$
with generators again denoted  $(e_+,h,e_-,e_+^*,h^*,e_-^*)$, and relations
$$
[h,e_\pm] = \pm 2e_\pm, [h^*,e_\pm^*] = 2e_\pm^*, [h, e_\pm^*] = \pm
2e_\pm^* , [h^*,e_\pm] = -2(e_\pm - 2e_\pm^*),
$$
$$
[e_+,e_-] = {{q^h - q^{-h}}\over{q - q^{-1}}}, \
[e_+,e_-^*] = {{q^{h} - q^{h^*}}\over{q - q^{-1}}}, \ 
[e_+^*,e_-] = {{q^{h^*} - q^{-h}}\over{q - q^{-1}}}, \ 
[e_+^*,e_-^*] = 0, 
$$
$$
[h,h^*]= 0, \quad e_\pm^* e_\pm - q^{-2} e_\pm e_\pm^* = (1-q^{-2})(e^*_\pm)^2 , 
$$
where we set $q = \exp(\hbar)$. 
We define $A_+$ (resp., $A_-$) to be the subalgebra of $A$ generated by $e_+,h,e_-$ 
(resp., $e_+^*,h^*,e_-^*$). There is a unique algebra map $\Delta : A \to A\otimes A$, 
such that 
$$
\Delta(e_+) = e_+\otimes q^{h} + 1\otimes e_+, \Delta(e_-) = e_-\otimes 1 + q^{-h}\otimes e_-,
\Delta(h) = h\otimes 1 + 1\otimes h, 
$$ 
and 
$$
\Delta(e_+^*) = (e_+^*\otimes q^{h^*}) \left( 
1 - q^{-1}(q - q^{-1})^2 e_+^*\otimes e_-^* \right)^{-1}
+ 1 \otimes e_+^*,  
$$
$$
\Delta(e_-^*) = e_-^* \otimes 1 + (q^{h^*}\otimes e_-^*)
\left( 1 - q^{-1}(q - q^{-1})^2 e_+^*\otimes e_-^* \right)^{-1}, 
$$
$$
\Delta(q^{h^*}) = (q^{h^*}\otimes q^{h^*})
\left( 1 - q^{-1}(q - q^{-1})^2 e_+^*\otimes e_-^* \right)^{-1}
\left(1 - q^{-3}(q - q^{-1})^2 e_+^*\otimes e_-^*\right)^{-1} . 
$$
 
Then $A_+$ and $A_-$ are Hopf subalgebras of $A$, and
$(A,A_+,A_-,\Delta)$ is a quantization of the Manin triple
$(\G,\G_+,\G_-)$. The algebra $A$ is the Drinfeld double of
$(A_+,\Delta)$ and its $R$-matrix is
$$
\cR = \exp_{q^2}\left( - (q - q^{-1}) e_-^* \otimes e_+\right) 
q^{{1\over 2} h^* \otimes h}
\exp_{q^2}\left( - (q - q^{-1}) e_+^* \otimes e_-\right)^{-1},  
$$ 
where $\exp_{q^2}(z) = \sum_{n\geq 0} {{z^n}\over{[n]!}}$ and $[n]! =
\prod_{k=1}^n (1 + q^2 + \cdots + q^{2k - 2})$. 

Let us fix $\al\in\CC$ and define $\HH_\al$ to be the subalgebra 
$\Ad(e^{\al h^*})(\G_+)$  of $\G$. The linear space $\HH_\al$ is 
spanned  by $h$ and  $e_{\pm} + \beta  e^*_{\pm}$, where $\beta =
e^{4\al} - 1$.  When  $\beta\neq 0$,  $(\G,\G_\pm,\G_\mp,\HH_\al)$ are 
nongraded Manin  quadruples. We now assume $\beta\neq 0$.  

Here are quantizations of these Manin quadruples. Let us define  $B_\al$
to be the subalgebra of $A$ generated by $h$ and $e_{\pm} + \beta 
e_{\pm}^*$. Let us set 
$$
F_\al = \Psi_\al(e_+^*\otimes e_-^*), \ \on{with}\ \Psi_\al(z) = 
{{\exp_{q^2} \left( -(q-q^{-1})e^{-4\al}z \right)
}\over{ \exp_{q^2} \left( -(q-q^{-1})z \right) }} .  
$$  

\begin{prop}
$((A,\Delta,\cR),A_\pm,B_\al,F_\al)$ are quantizations of the Manin quadruples
$(\G,\G_\pm,\G_\mp,\HH_\al)$. 
\end{prop}

{\em Proof.}  We have $B_\al\cap A_+ = \CC[h][[\hbar]]$, and $B_\al\cap
A_- = \CC[[\hbar]]$. 

Let us set $u_\al = e^{\al h^*}$.  Then $\Delta(u_\al) = F^*
(u_\al\otimes u_\al) (F^*)^{-1}$,  where $F^* = \exp_{q^2} 
(-(q-q^{-1})
e_+^*\otimes e_-^*)$. It follows that  $F_\al = (u_\al\otimes
u_\al) \Delta(u_\al)^{-1}$, which  implies that $F_\al$ satisfies the
cocycle identity. 

Moreover, $B_\al = u_\al A_+ u_\al^{-1}$, therefore 
$F_\al\Delta(B_\al)F_\al^{-1}\subset B_\al^{\otimes 2}$. In fact, $A$ 
equipped with the twisted  coproduct $F_\al\Delta F_\al^{-1}$ is a Hopf
algebra, and thus $B_\al$ is a Hopf subalgebra of $(A,F_\al \Delta
F_\al^{-1})$. 

Let us now show that $F_\al$ satisfies both condition (\ref{hyp}), and
the similar condition where $A_+$ is replaced by $A_-$. 
This follows from the conjunction of 
\begin{equation} \label{cond:1}
\Delta(u_\al) \in  
\left( 1 + \hbar (A_-)_0\otimes (A_+)_0 \right) (u_\al\otimes u_\al)
\left( 1 + \hbar A_0 \otimes (A_+)_0 \right)  
\end{equation}
and
\begin{equation} \label{cond:2}
\Delta(u_\al)\in \left( 1 + \hbar (A_+)_0\otimes (A_-)_0 \right) 
(u_\al\otimes u_\al)
\left( 1 + \hbar (A_+)_0 \otimes A_0 \right) . 
\end{equation} 
Then 
\begin{equation} \label{eq:twists}
\Delta(u_\al) = F' (u_\al\otimes u_\al) (F')^{-1}
= F'' (u_\al\otimes u_\al) (F'')^{-1} , 
\end{equation}
where 
$$
F' = \exp_{q^2} ( - (q - q^{-1}) e_+^*\otimes e_-), \quad 
F'' = \exp_{q^2} ( - (q - q^{-1}) e_+\otimes e_-^*).    
$$
The first equality of (\ref{eq:twists}) proves (\ref{cond:1}), and the
second one proves (\ref{cond:2}).  \hfill \qed\medskip

One may expect that the algebra-module over $(A_+^{opp},\Delta)$
constructed by means of $((A,\Delta,\cR),A_+,B_\al,F_\al)$ in Theorem
\ref{thm:bis} is a formal completion of the function algebra over a
Podle\`s sphere (\cite{Podles}).

\begin{remark} 
One shows that for $x\in A_-$, $\Delta(x) = F^* \wt\Delta(x)
(F^*)^{-1}$, where   $\wt\Delta$ is the coproduct on $A_-$ defined by
$\wt\Delta(e_+^*) = e_+^*\otimes q^{h^*} + 1 \otimes e_+^*$, 
$\wt\Delta(q^{h^*}) = q^{h^*} \otimes q^{h^*}$, and  $\wt\Delta(e_-^*) =
e_-^*\otimes 1 + q^{h^*} \otimes e_-^*$.  The completion of 
the Hopf algebra $(A_-,\Delta)$ with respect to the topology 
defined by its augmentation ideal 
should be isomorphic to the formal completion at the identity of the quantum
coordinate ring of $SL_2$.  The Hopf algebra $A$ was obtained  by a
method similar to that of Drinfeld's ``new realizations''.  
\end{remark}

\subsection{The case of loop algebras}

Let us return to the situation of the Manin quadruple 
(\ref{infinite:quadruple}), in the case where $\cK$ is a field of 
Laurent series and $R$ is a ring of functions on an affine  curve. In
this situation, one can consider the following problems. 

1) If the double quotient $G_+\setminus G / H$ is equipped with the zero
Poisson structure, the projection $G/H \to G_+ \setminus G / H$ is
Poisson.  The ring of formal functions on this double quotient is  
$\cO_{G_+ \setminus G / H} = \left( U\G / ((U\G) \HH + \G_+ U\G)
\right)^*$. On the other hand, $(A^*)^{A_+,B} = \{\ell\in A^* | \forall
a_+\in A_+, b\in B, \ell(a_+ a b) = \ell(a)\eps(a_+)\eps(b)\}$ is a
subalgebra of $(A^*)^B$.  It is commutative because the $R$-matrix
$\cR$ of $(A,\Delta)$ belongs to $(A_+\otimes A_-)^\times_0$ and the
twisted $R$-matrix $F^{(21)}\cR F^{-1}$ belongs to $(B\otimes
A)^\times_0$ (see \cite{ER}).

It would be interesting to describe the algebra inclusion $(A^*)^{A_+,B}
\subset (A^*)^B$, to see whether it is a flat deformation of $\cO_{G_+
\setminus G / H}$,  and when the level is critical, to describe the
action of the quantum Sugawara field by commuting operators on
$(A^*)^{A_+,B}$ and $(A^*)^B$. These operators could be related to the
operators constructed in \cite{EF}. 

2) For $g$ fixed in $G$, let $\null^g x$  denote the conjugate
$gxg^{-1}$ of an element $x$ in $G$ by $g$.  One would like to describe
the Poisson  homogeneous spaces $G_\pm / (G_\pm \cap \null^g H)$, and to
obtain  quantizations of the Manin quadruples
$(\G,\G_\pm,\G_\mp,\null^g\HH)$.  A natural idea would be to start from
the quantization of the Manin  quadruple   $(\G,\G_\pm,\G_\mp,\HH)$
obtained in \cite{ER}, and to apply to $B$  a suitable automorphism of
$A$ which lifts the automorphism $\Ad(g)$ of $U\G$. 

We hope to return to these questions elsewhere.

\end{document}